\documentclass[english,11pt,a4paper,twoside]{article}

\usepackage[english]{babel}
\usepackage{enumitem} % ermoeglicht manuelle Aenderung von itemize und enumerate
\usepackage{wrapfig}
\usepackage[final]{graphicx}
\usepackage{epsfig}
\usepackage{textcomp}
\usepackage{amsmath}
\usepackage{amssymb}
\usepackage{amsfonts}
\usepackage{mathbbol} % fuer Befehl \Eins
\usepackage{psfrag}
\usepackage{latexsym}
\usepackage{eurosym}
\usepackage{color}    % \definecolor,\textcolor
\usepackage{mathtools}
\usepackage{tikz}
\usepackage{stmaryrd} % fuer dreifache Normstriche mit \interleave
\usepackage[normalem]{ulem}
\usepackage{verbatim} % fuer comment Befehl

\usepackage{hyperref}
\hypersetup{colorlinks=true}
\hypersetup{linkcolor=blue}
\hypersetup{citecolor=blue}

\newcommand{\epsi}{\varepsilon}
\newcommand{\IC}{\mathbb{C}} 

\newcommand{\IN}{\mathbb{N}}

\newcommand{\IR}{\mathbb{R}}

\newcommand{\tend}{\tt_{\mbox{\tiny end}}}

\newcommand{\ds}{\displaystyle}
\newcommand{\hs}[1]{\hspace*{#1 mm}}

\newcommand{\Ord}[1]{\mathcal{O} \! \left(#1\right)}

\newcommand{\A}{\mathcal{A}}
\newcommand{\F}{\mathcal{F}}

\newcommand{\J}{\mathcal{J}}
\renewcommand{\L}{\mathcal{L}}
\renewcommand{\P}{\mathcal{P}}
\newcommand{\Pp}{\mathcal{P}^\perp}
\renewcommand{\S}{\mathcal{S}}

\newcommand{\ii}{\mathrm{i}}

\newcommand{\w}{\omega}

\newcommand{\uu}{u}

\newcommand{\vv}{v}
\newcommand{\vvtilde}{\widetilde{v}}
\newcommand{\cg}{c_{g}}

\newcommand{\pt}{\partial_{\tt}}

\newcommand{\intl}{\int\limits}

\newcommand{\stkout}[1]{\ifmmode\text{\sout{\ensuremath{#1}}}\else\sout{#1}\fi}

\newcommand{\weg}[1]{}
\newcommand{\proof}{\noindent\textbf{Proof. }}

\def\qed {\hfill \rule{0,2cm}{0,2cm} \bigskip \\}

\newcommand{\eps}{\varepsilon}
\def\t{t}
\def\tt{\tau}
\newcommand{\x}{x}
\newcommand{\xx}{\xi}
\newcommand{\pp}{p}
\newcommand{\RR}{R}
\newcommand{\err}{e}

\newtheorem{Theorem}{Theorem}[section]

\newtheorem{Lemma}[Theorem]{Lemma}
\newtheorem{Corollary}[Theorem]{Corollary}
\newtheorem{Assumption}[Theorem]{Assumption}

\newtheorem{Remark}[Theorem]{Remark}

\newcounter{fig}

\begin{document}

\title{Polarized high-frequency wave propagation beyond the nonlinear Schr\"odinger approximation\thanks{Funded by the Deutsche Forschungsgemeinschaft (DFG, German Research Foundation) – Project-ID 258734477 – SFB 1173}}
\author{Julian Baumstark\thanks{
Karlsruher Institut f\"ur Technologie, Fakult\"at f\"ur Mathematik, Institut für Angewandte und Numerische Mathematik, \texttt{julian.baumstark@kit.edu}, \texttt{tobias.jahnke@kit.edu}
} \and Tobias Jahnke\footnotemark[2] \and Christian Lubich\thanks{Universit\"at T\"ubingen, Mathematisches Institut, \texttt{lubich@na.uni-tuebingen.de}}
}
\date{\today}
\maketitle

\begin{abstract}
This paper studies highly oscillatory solutions to a class of systems of semilinear hyperbolic equations with a small parameter, in a setting that includes Klein--Gordon equations and the Maxwell--Lorentz system. The interest here is in solutions that are polarized in the sense that up to a small error, the oscillations in the solution depend on only one of the frequencies that satisfy the dispersion relation with a given wave vector appearing in the initial wave packet. 
The construction and analysis of such polarized solutions is done using modulated Fourier expansions. This approach includes higher harmonics and yields approximations to polarized solutions that are of arbitrary order in the small parameter, going well beyond the known first-order approximation via a nonlinear Schr\"odinger equation. The given construction of polarized solutions is 
explicit, uses in addition a linear Schr\"odinger equation for each further order of approximation, and is accessible to direct numerical approximation.
\end{abstract}

% --------------------------------------------------------------------------------
\section{Introduction}
% --------------------------------------------------------------------------------

We consider semilinear hyperbolic systems of the form
\begin{align}
\label{PDE.uu}
 \partial_\t  \uu +A(\partial_\x)\uu +\frac{1}{\epsi} E\uu &= \epsi \, T(\uu,\uu,\uu),
 \qquad \x \in\IR^d, \,\t \in [0,\tend/\epsi]
\end{align}
with vector-valued solution $\uu(\t ,\x ) \in \IR^n$ ($d,n \in \IN$).
Here, $ A(\partial_\x) = \sum_{\ell=1}^d A_\ell \partial_\ell $ is a first-order differential operator with symmetric matrices $A_1, \ldots, A_d \in \IR^{n \times n}$,  
$ E\in\IR^{n\times n} $ is a skew-symmetric matrix,
$T: \IR^n \times \IR^n \times \IR^n \rightarrow \IR^n $ is a trilinear nonlinearity, and $0<\eps\ll 1$ is a small parameter. 
%Its trilinear extension to $\IC^n \times \IC^n \times \IC^n$, which will be considered below, is also denoted by $T$.
This problem class includes, e.g., the Maxwell--Lorentz system and Klein--Gordon systems; 
cf.~Section 2.1 in \cite{colin-lannes:09}.

In \cite{colin-lannes:09,lannes:11} such systems have been considered with initial data
\begin{align}
\label{PDE.uu.c}
 \uu(0,\x ) &= e^{\ii \kappa \cdot \x /\eps} p(\x ) + c.c.,
\end{align}
where $\kappa \in \IR^d\setminus\{0\}$ is a fixed wave vector and $p: \IR^d \rightarrow \IR^n$ a smooth envelope independent of $\eps$. (Different scalings such as $p(x/\sqrt\eps)$ instead of $p(x)$ are also of interest, but will not be considered here.) The dot $\cdot$ denotes the Euclidean scalar product in $\IR^d$ and ``c.c.'' means complex conjugation of the preceding term.

Both the PDE \eqref{PDE.uu} and the initial data \eqref{PDE.uu.c} contain the small parameter
$0<\epsi\ll 1$, and as a consequence, solutions of \eqref{PDE.uu}--\eqref{PDE.uu.c} typically
oscillate rapidly in time and space. 

\paragraph{Polarization in the linear case.}
When $T\equiv 0$ in \eqref{PDE.uu}, a situation of particular interest appears for special initial data $p(x)=\alpha(x)v$, where $\alpha$ is a smooth real- or complex-valued function, and $v$ is an eigenvector to an eigenvalue\footnote{In the introduction we assume $\w$ to be a simple eigenvalue.} $\omega$ of the Hermitian matrix $\sum_{\ell=1}^d \kappa_\ell A_\ell - \ii E$ (in different terminology, the pair of the wave vector $\kappa$ and the frequency $\w$ satisfies the dispersion relation).  Then, the solution is of the form
$u(t,x)=e^{\ii (\kappa \cdot \x  - \w \t )/\epsi} \mu_{L}(t,x) v + c.c.$, where $\mu_{L}$ is a smooth, $\eps$-independent complex-valued function. This is a \emph{polarized solution} where the other eigenvalues $\widetilde\omega \ne \omega$ do not appear in highly oscillatory exponentials
$e^{\ii (\kappa \cdot \x  - \widetilde\w \t )/\epsi}$, as they would for general initial data.

\paragraph{Nonlinear polarization.}
In the present paper we contribute to extending the concept of polarization to the \emph{nonlinear} situation \eqref{PDE.uu}. 
There is earlier work by Colin \& Lannes~\cite{colin-lannes:09} and Lannes~\cite{lannes:11}, which motivated the present work. They considered initial data \eqref{PDE.uu.c} with $p(x)=\alpha(x)v + \Ord{\eps}$ and  showed that the solution can then be written as
$u(t,x)=e^{\ii (\kappa \cdot \x  - \w \t )/\epsi} \mu(t,x) v + c.c. + \Ord{\eps}$
with a function $\mu$ that, together with its partial derivatives, is bounded independently of $\eps$. (Here and in the following, the $\Ord{\cdot}$ notation is with respect to the maximum norm.) This approximation result holds true over times $\Ord{\eps^{-1}}$, which is the time scale over which the nonlinearity yields an $O(1)$ contribution to the solution.

Here, we are interested in approximations of higher order in $\eps$. We will
show that for given smooth $p(x)=\alpha(x)v$, there exists a solution of \eqref{PDE.uu} of the form (for an arbitrary fixed odd integer $m\ge 1$)
\begin{equation}
\label{mfe-intro}
\uu(\t,\x) = \sum_{|j|\le m \atop j \text{ odd}} \uu_j(\t,\x) \, e^{\ii j (\kappa\cdot \x - \omega \t)/\eps} + \Ord{\eps^{m+1}}
\end{equation}
such that $\P\uu_1(0,\x)=p(x)$ with the orthogonal projection $\P=vv^*$ onto the eigenspace, and the modulation functions $\uu_j$ with $\uu_{-j}=\overline \uu_j$ are smooth and,  together with their partial derivatives of any order, they are bounded by  $\uu_j=\Ord{\eps^{|j|-1}}$. The dominant modulation function $u_1(t,x)$ is aligned with the eigenvector $v$ up to $\Ord{\eps}$. Moreover, the representation \eqref{mfe-intro}
 is unique up to $\Ord{\eps^{m+1}}$. This higher-order approximation result again holds true over times $\Ord{\eps^{-1}}$.
 Note that here no other frequencies $\widetilde\w \ne \w$ are allowed to appear in the oscillatory exponentials in \eqref{mfe-intro}. It is not obvious {\it a priori} that this can be achieved for the semilinear equation.
 
\paragraph{Nonlinear Schr\"odinger approximation.}
In
Section 2.3.2 of \cite{colin-lannes:09} and Section 2.3 of \cite{lannes:11},
the solution of \eqref{PDE.uu}--\eqref{PDE.uu.c} with $p(x)=\alpha(x)v$ is approximated up to $\Ord{\eps}$ by
\begin{align}
\label{u.NLS}
\uu_{\text{\tiny NLS}}(\t ,\x ) = 
e^{\ii (\kappa \cdot \x  - \w \t )/\epsi} U(\eps\t ,\x -c_g\t) + c.c.,
\end{align}
where $ U(\tt,\xx) $ is the solution of an $\eps$-independent nonlinear Schr\"odinger equation (NLS) on the rescaled time interval $[0,\tend]$. In the second argument of~$U$, $\cg$ is the group velocity.

The modulation function $u_1$ in \eqref{mfe-intro} admits an asymptotic expansion in powers of $\eps$: \ 
$u_1=u_1^0+\eps u_1^1+\eps^2 u_1^2 + \dots$\ .
It turns out that $u_1^0(t,x)$ in \eqref{mfe-intro} equals $U(\eps\t ,\x -c_g\t) $, though it arises with a different proof.

\paragraph{Higher harmonics.}
A shortcoming of the nonlinear Schr\"odinger approximation is the fact that substituting \eqref{u.NLS} into the nonlinearity yields
\begin{align*}
T(\uu_{\text{\tiny NLS}},\uu_{\text{\tiny NLS}},\uu_{\text{\tiny NLS}})
=&\
e^{\ii (\kappa \cdot \x  - \w \t )/\epsi} 
\Big(T(U,U,\overline{U})+T(U,\overline{U},U)+T(\overline{U},U,U)\Big)
\\
&\!\!\!+
e^{-\ii (\kappa \cdot \x  - \w \t )/\epsi} 
\Big(T(U,\overline{U},\overline{U})+T(\overline{U},U,\overline{U})+T(\overline{U},\overline{U},U)\Big)
\\
&\!\!\!+ e^{3\ii (\kappa \cdot \x  - \w \t )/\epsi} T(U,U,U)
+
e^{-3\ii (\kappa \cdot \x  - \w \t )/\epsi} T(\overline{U},\overline{U},\overline{U}),
\end{align*}
whereas substituting \eqref{u.NLS} into the left-hand side of \eqref{PDE.uu} does not produce any terms with factor $e^{\pm 3\ii (\kappa \cdot \x  - \w \t )/\epsi}$. Hence, these \emph{higher harmonics} contribute to the residual and thus to the approximation error. In order to improve the accuracy, these higher harmonics have to be taken into account. This is a first motivation to search for approximations to $\uu$ of the form \eqref{mfe-intro}.

\paragraph{Modulated Fourier expansions and geometric optics.}
Asymptotic expansions such as \eqref{mfe-intro} are known as modulated Fourier expansions (MFEs) in time, a name that was coined in
\cite[Chapter XIII]{HaiLW02}. Multiscale expansions of this type were occasionally used before under different names: in a formal way without rigorous bounds in \cite{Kru58} for analysing the motion of charged particles and in
\cite{Whi74} for asymptotics of wave propagation, and with rigorous bounds of modulation functions and remainder term in \cite{HaiL00} for studying the long-time behaviour of numerical methods for highly oscillatory Hamiltonian systems. Since then, MFEs were used in many works on the analysis of highly oscillatory ordinary and evolutionary partial differential equations and their numerical approximations; see, e.g., the reviews \cite{HaiL12,GauHL18} and numerous references therein. 

In the way MFEs are used here, with an expansion of the modulation functions in powers of $\eps$ (which is possible for the problem at hand but not in other cases with near-resonances between frequencies; see e.g. \cite{CohHL08,FaoGL13}), the approach of this paper is actually close to the often-used approach in geometric optics of making an ansatz (see e.g. \cite[Chapters~7--9]{rauch-book:12}) 
\begin{align*}
u(t,x) = \sum_{j} \epsi^j a_j(t,\eps t,x,(\kappa \cdot \x  - \widetilde\w \t )/\epsi),
\end{align*}
where $ \theta \mapsto a_j(t,\tau,x,\theta) $ is $2\pi$-periodic. Expanding this map into a Fourier series, as is often done, 
and interchanging the appearing double sum yields an expression that is formally very similar 
to a modulated Fourier expansion \eqref{mfe-intro} with modulation functions $u_j(t,x)$ expanded in powers of $\eps$. The actual construction of the coefficient functions and the derivation of bounds for them uses, however, different perspectives and details. 

For systems similar to or more general than \eqref{PDE.uu} asymptotic expansions of solutions have been constructed in 
\cite{rauch-book:12,donnat-rauch:97b,joly-metivier-rauch:00,joly-metivier-rauch:93} and many other works.
These articles refer to the setting of geometric optics, where the time interval does not depend on $\epsi$. This is to be distinguished from diffractive geometric optics, where the length of the time interval is
$\Ord{1/\epsi}$ as in \eqref{PDE.uu}.
In the diffractive regime, solutions of the initial-value problem for semilinear and quasilinear systems with $\eps$-oscillatory initial data \eqref{PDE.uu.c} have been approximated with arbitarily high order in $\eps$ in \cite{donnat-joly-metivier-rauch:96}, though with $\epsi E$ instead of $E/\epsi$ in \eqref{PDE.uu}. 
For \eqref{PDE.uu} with initial data \eqref{PDE.uu.c} a second-order approximation of the form
\eqref{mfe-intro}, but with oscillating $u_j$, has been analyzed in \cite{BauJah:22}.
Note that the problem of nonlinear polarization posed in the present paper is different from the standard initial value problem.

\paragraph{Computational aspects.}
Because of the high-frequency oscillations, a direct discretization of \eqref{PDE.uu}--\eqref{PDE.uu.c} would have to use very small time steps and spatial mesh widths and is therefore computationally prohibitive, and even more so the computation of polarized solutions.
The smooth, non-oscillatory modulation functions $u_j$ in \eqref{mfe-intro} are constructed explicitly in this paper, requiring only the solution of one nonlinear and a few linear Schr\"odinger equations that are independent of $\eps$, using co-moving coordinates and the slow time scale $\tau=\eps t$, on the $\eps$-independent rescaled time interval $[0,\tend]$. 
With the NLS approximation or the MFE, one can thus achieve 
an accuracy of $\Ord{\epsi}$ or $\Ord{\eps^{m+1}}$, respectively, without having to solve any oscillatory PDEs. While computational aspects were the basic motivation for the present paper, they will not be elaborated here.

\paragraph{Outline.} In Section~\ref{Sec.result} we give a precise formulation of the setting and present our main result, Theorem~\ref{thm.mfe}, which establishes the polarized modulated Fourier expansion \eqref{mfe-intro}.
For $m=1$ equations for the computation of the modulation functions are stated explicitly in Corollary~\ref{Corollary.example.m=1}.
The remainder of the paper gives the proof of Theorem~\ref{thm.mfe}. In Section~\ref{Sec.construction} we construct the MFE formally, in Section~\ref{Sec.bounds} we derive bounds for the modulation functions, and in Section~\ref{Sec.error} we bound the approximation error of the truncated MFE.

% --------------------------------------------------------------------------------
\section{Statement of the main result}
\label{Sec.result}
% --------------------------------------------------------------------------------

For every $k \in \IR^d$ we define the real symmetric matrix
\begin{align}
\label{Def.Ak}
 A(k)&=\sum_{\ell=1}^d k_\ell A_\ell \in \IR^{n\times n}
\end{align}
with the same matrices $ A_\ell$ as in the differential operator $A(\partial_\x)$.
A pair $(\w,k)$ fulfills the \emph{dispersion relation} if $\ii\w/\eps$ is an eigenvalue 
of $A(\ii k/\eps)+E/\eps$, or equivalently, $\w$ is an eigenvalue of the Hermitian matrix $A(k)-\ii E$.

\begin{Assumption} \label{Ass.om}
Let $\kappa \in \IR^d\setminus\{0\}$ be a fixed wave vector, and let 
$\w=\w_1(\kappa)$ be a selected eigenvalue of $A(\kappa)-\ii E$. We assume that $\w_1(\cdot)$ is twice continuously differentiable 
at~$\kappa$.
\end{Assumption}

We note that a sufficient (though not necessary) condition for the required differentiability is that $\w$ is a simple eigenvalue of $A(\kappa)-\ii E$.

\medskip
The Hermitian matrix
\begin{align*}
\L(\w,\kappa) =  -\w I + A(\kappa)-\ii E \in \IC^{n\times n}
\end{align*}
thus has a nontrivial kernel, which is the eigenspace to the eigenvalue $\w$ of $A(\kappa)-\ii E$.
Let $\P$ be the corresponding eigenprojection, i.e., the orthogonal projection onto the kernel of $\L(\w,\kappa)$. We set $\Pp=I-\P$.

\medskip

We consider a function $\pp$ that takes values in the range of $\P$:
\begin{equation}\label{u-init-par}
\text{$\pp\in \S(\IR^d,\IC^n)$ \ with \  $\P \pp(\x)=\pp(\x)$  for all  $\x\in\IR^d$.}
\end{equation}
Here, $\S(\IR^d,\IC^n)$ is the Schwartz space of smooth functions that, together with all their derivatives, decay faster than every negative power of $|x|$. The function $\pp$ is assumed to be independent of $\eps$.

\medskip
Our next assumption concerns the NLS approximation of \eqref{u.NLS}, $\uu_{\text{\tiny NLS}}(\t ,\x ) = 
e^{\ii (\kappa \cdot \x  - \w \t )/\eps} U(\eps\t ,\x - c_g \t) + c.c.$  with initial data $\uu_{\text{\tiny NLS}}(0,\x)= \pp(\x) e^{\ii \kappa\cdot\x} + c.c.$, where $c_g=\nabla\w_1(\kappa)$ is the group velocity.
With the Hessian of $\w_1$ at $\kappa$, denoted $H=\nabla^2\w_1(\kappa)$, the function $U=U(\tt,\xx)$ is taken as the solution of the NLS
\begin{equation}
\label{U.NLS}
\begin{aligned}
& \partial_\tt U
 -\frac{\ii}{2}\nabla_\xx \cdot H\nabla_\xx U  
 = \P\big(T(\overline U,U,U)+ T(U,\overline U,U)+T(U,U,\overline U) \big)
 \\
& \text{with initial data $U(0,\cdot)=p$ of \eqref{u-init-par}.}
\end{aligned}
\end{equation}
Note that $\eps$ does not show up in \eqref{U.NLS}.
We assume the following. 

\begin{Assumption}\label{Ass.NLS} The NLS \eqref{U.NLS} has a bounded solution $U$ on $[0,\tend]\times\IR^d$ that
has arbitrarily many bounded partial derivatives.
\end{Assumption}

This is satisfied at least for sufficiently small $\tend>0$, as follows from the local existence result in \cite{Seg63}.

\medskip
For the construction of higher-order approximations in powers of $\varepsilon$ we need the following assumption.

\begin{Assumption}\label{Ass.Lj.invertible}
The matrix
\begin{align}
\label{Def.Lj}
\L_j = \L(j\w,j\kappa) \in\IC^{n\times n}
\end{align}
is invertible for every odd integer $j$ with $3 \leq j \leq m+2$.
\end{Assumption}
For $j=3$ this assumption was also made in \cite[Assumption 3]{colin-lannes:09}.

\medskip
Under these assumptions we have the following theorem, which is the main result of this paper.

\begin{Theorem}[Modulated Fourier expansion] \label{thm.mfe}
Let $\kappa \in \IR^d\setminus\{0\}$ be a fixed wave vector and let the frequency $\omega=\omega_1(\kappa)\in\IR$. 
Let $\pp\in \S(\IR^d,\IC^n)$ with $\P \pp=\pp$. Fix an odd positive number $m$.
Under Assumptions~\ref{Ass.om}--\ref{Ass.Lj.invertible}, there exists a solution
$\uu$ of \eqref{PDE.uu} on the interval $[0,\tend/\varepsilon]$ that admits a modulated Fourier expansion
\begin{equation}\label{mfe}
\uu(\t,\x) = \sum_{|j|\le m \atop j \text{ odd}} \uu_j(\t,\x) \, e^{\ii j (\kappa\cdot \x - \omega \t)/\eps} + \RR_{m}(\t,\x)
\end{equation} 
with $\uu_{-j}(\t,\x)=\overline{\uu_j(\t,\x)}$ and with the following properties:
\begin{enumerate}
\item
\label{thm.mfe.property.Puu}
$\P \uu_1$ satisfies the initial condition
$$
\P \uu_1(0,\cdot)=\pp.
$$
\item The modulation functions $\uu_j$  are bounded as follows, with a constant $C$ independent of~$\eps$ (but depending on~$m$):
\begin{align}
\label{u1-bounds}
&\| \P \uu_1 \|_{L^\infty([0,\tend/\varepsilon]\times \IR^d)} \le C, \qquad
\|  \P^\perp \uu_1 \|_{L^\infty([0,\tend/\varepsilon]\times \IR^d)} \le C\eps,
\\
\label{uj-bounds}
& \|  \uu_j \|_{L^\infty([0,\tend/\varepsilon]\times \IR^d)} \le C\eps^{j-1} \qquad (\text{for odd $\,3 \le j\le m$}),
\end{align}
and the same bounds hold true for the spatio-temporal partial derivatives of these functions up to any fixed order.
Moreover, the modulation functions $\uu_j$ admit an asymptotic expansion in powers of $\eps$. 

\item The remainder term is bounded by
\begin{equation}
\label{R-bound}
 \| \RR_{m} \|_{L^\infty([0,\tend/\varepsilon]\times \IR^d)} \le C \eps^{m+1}.
\end{equation}
\end{enumerate}
The modulation functions $\uu_j$ with these properties (and hence also $\uu$) are unique up to $O(\eps^{m+1})$ w.r.t. the maximum norm.
\end{Theorem}

\begin{Remark}[Polarization]
We note in particular that the initial value 
$ \uu(0,\x) = (e^{\ii \kappa \cdot \x /\eps} p(\x ) + c.c.) + O(\eps) $ is determined by $\pp$ uniquely up to $O(\eps^{m+1})$. So Theorem~\ref{thm.mfe} holds true for special initial values $\uu(0,\cdot)$ of \eqref{PDE.uu}, which we call \emph{polarized initial values}. As will be seen from the explicit construction in the next section, the polarization map $\pp \mapsto \uu(0,\cdot)$ is a nonlinear map.  The construction shows that polarized initial values correspond to a slow manifold of a system of highly oscillatory partial differential equations.

\smallskip
For polarized initial values, the other frequencies $\omega_2(\kappa),\dots,\omega_\mu(\kappa)$ 
(i.e.\ the other eigenvalues of $A(\kappa)-\ii E$)
do not enter the oscillatory exponentials of the modulated Fourier expansion. For general initial values of type \eqref{PDE.uu.c}, it would be possible to obtain a modulated Fourier expansion that includes also those frequencies, under appropriate resonance or non-resonance conditions among the frequencies. However, this is not the point of the above result, where the interest lies precisely in the fact that the other frequencies do \emph{not} occur.
\end{Remark}

\begin{Corollary}\label{Corollary.example.m=1}
Choosing $m=1$ in \eqref{mfe} gives an $\Ord{\eps^2}$ approximation on the time interval $[0,\tend/\eps]$ that is constructed as follows:
We obtain $u_1= p_1 +q_1$ (with $\P u_1= p_1$ and $\Pp u_1=q_1$) by setting 
$$
p_1=p_1^0+\eps p_1^1, \qquad q_1 =\eps q_1^1,
$$
where 
\begin{itemize}
\item $p_1^0(t,x) = U(\eps t, x- \cg t)$, where $U(\tau,\xi)$ is the solution to the nonlinear Schr\"odinger equation \eqref{U.NLS} with initial value $U(0,\cdot)=p$.
\item $q_1^1(t,x) = z_1^1(\eps t, x- \cg t)$, where $z_1^1(\tau,\xi)$ is obtained from  $\partial_\xi U(\tau,\xi)$ via matrix operations detailed in \eqref{z11}.
\item $p_1^1(t,x) = y_1^1(\eps t, x- \cg t)$, where $y_1^1(\tau,\xi)$ is the solution to the linear Schr\"odinger equation \eqref{y11.LS} with zero initial value.
\end{itemize}
\end{Corollary}

\begin{Remark}[Computation of polarized solutions]
Since the proof of Theorem~\ref{thm.mfe} gives an explicit construction of the modulation functions, the result can be used to devise numerical schemes to compute polarized solutions of \eqref{PDE.uu}, using time steps that are much larger than $\eps$. This will be presented elsewhere.
\end{Remark}

% --------------------------------------------------------------------------------
\section{Construction of the modulated Fourier expansion}
% --------------------------------------------------------------------------------
\label{Sec.construction}

% --------------------------------------------------------------------------------
\subsection{Rescaling of time and co-moving coordinate system}
% --------------------------------------------------------------------------------

For our analysis it is convenient to change to the new variables
\begin{align}
\label{Comoving.coordinates}
\tt=\epsi \t , \qquad \xx=\x -\cg \t , \qquad
 \vv(\tt,\xx) = \vv(\epsi\t ,\x -\cg \t) = \uu(\t ,\x ),
\end{align}
where again $\cg = \nabla \w_1(\kappa)$ is the group velocity. Since
\begin{align*}
\partial_\t  \uu(\t ,\x ) 
&= 
\partial_\t  \vv(\epsi\t ,\x -\t \cg)
=
\epsi \partial_\tt \vv(\tt,\xx)
- \cg \cdot \nabla_\xx \vv(\tt,\xx), 
\\
A(\partial_\x )\uu(\t ,\x ) 
&= A(\partial_\xx)\vv(\tt,\xx),
\end{align*}
the change of variables turns the original problem \eqref{PDE.uu} into
\begin{align}
\label{PDE.vv}
 \pt \vv + \frac{1}{\epsi}B(\partial_\xx)\vv +\frac{1}{\epsi^2} E\vv &= T(\vv,\vv,\vv),
 \qquad \xx\in\IR^d, \; \tt\in [0,\tend]
\end{align}
with 
\begin{align}
\label{Def.B}
 B(\partial_\xx) = \sum_{\ell=1}^d B_\ell \frac{\partial}{\partial \xx_\ell}, \qquad
 B_\ell = A_\ell - (\cg)_\ell I.
\end{align}
We choose a fixed odd number $m \geq 1$ and define the index sets
\begin{align}
\label{Def.J}
\J &= \{\pm 1, \pm 3, \ldots, \pm m \}\quad \text{and} \quad \J_+=\J \cap \IN.
\end{align}
We aim to construct modulation functions $ v_j(\tt,\xx) = u_j(\t ,\x ) $ such that in the new variables $(\tt,\xx)$,
\begin{align}
\label{Ansatz.v}
 \vv(\tt,\xx) &\approx
 \sum_{j\in \J}
e^{\ii j \kappa \cdot \xx/\epsi} 
e^{\ii j (\kappa \cdot c_g - \w) \tt/\epsi^2} 
v_j(\tt,\xx),
\qquad v_{-j} = \overline{v_j}
\end{align}
with the bounds in the maximum norm
\begin{subequations}
\label{Ass.initial.data}
\begin{align}
\label{Ass.initial.data.a}
\P v_1 &= \Ord{1}, 
\\
\label{Ass.initial.data.b}
\Pp v_1 &= \Ord{\epsi}, 
\\
\label{Ass.initial.data.c}
v_j &= \Ord{\epsi^{j-1}} 
\quad \text{for odd} \quad 3\leq j \le m,
\end{align}
\end{subequations}
and the same bounds for the space and time derivatives up to some fixed order, and with an  
\begin{equation} \label{v-err}
\text{approximation error in \eqref{Ansatz.v} of size $\Ord{\eps^{m+1}}$ w.r.t. maximum norm.}
\end{equation}
When this is achieved, this yields \eqref{mfe} in the original variables $(\t,\x)$ together with the bounds of Theorem~\ref{thm.mfe}.

% --------------------------------------------------------------------------------
\subsection{PDEs for $v_j$}
% --------------------------------------------------------------------------------

Substituting \eqref{Ansatz.v} into \eqref{PDE.vv}, collecting terms with the same exponential factor and discarding terms with factors
\begin{align*}
e^{\ii j \kappa \cdot \xx/\epsi} e^{\ii j (\kappa \cdot c_g - \w) \tt/\epsi^2},
\qquad |j|>m ,
\end{align*}
yields the system
\begin{align}
\label{PDE.vj}
 \pt v_j
 + \frac{\ii}{\epsi^2}\L_j v_j
 + \frac{1}{\epsi}B(\partial_\xx)v_j
 =
 \sum_{\#J=j} T(v_J),&
\\
\notag
j\in\J_+, \quad t\in[0,\tend], \quad x\in\IR^d &
 \end{align}
with $\L_j$ defined in \eqref{Def.Lj} and the notation $J=(j_1,j_2,j_3)$, 
$\#J:= j_1+j_2+j_3$,
and
$ T(v_J)=T(v_{j_1},v_{j_2},v_{j_3})$.
The sum on the right-hand side is taken over the set
 \begin{align*}
 \Big\{J=(j_1,j_2,j_3) \in \J^3 : \#J=j \Big\}.
\end{align*}
The system of PDEs \eqref{PDE.vj} is compatible with the condition that $ v_{-j} = \overline{v_j}. $
Since $ \L_j$ is invertible for $ 3 \leq j \in \J_+ $ according to Assumption~\ref{Ass.Lj.invertible},
\eqref{PDE.vj} can be reformulated as
\begin{align}
\label{PDE.vj.modified}
 v_j
 =
 \ii \epsi \L_j^{-1}
 \Big(\epsi \pt v_j + B(\partial_\xx)v_j - \epsi \sum_{\#J=j} T(v_J)\Big),
 \qquad 3 \leq j \in \J_+.
 \end{align}
The case $j=1$ is special.
We distinguish
\begin{subequations}
\label{Def.y1.z1}
 \begin{align}
 y_1 &= \P v_1, & y_{-1} &= \overline{y_1},
 \\
 z_1 &= \Pp v_1, & z_{-1} &= \overline{z_1}.
\end{align} 
\end{subequations}
In order to derive equations for $y_1$ and $z_1$, we use that $\P \L_1 = \L_1 \P = 0$ by definition, and that
$ \P A(\partial_\xx) \P = \P (c_g \cdot \nabla) $; cf.~\cite[Lemma 2.9]{lannes:11}. It follows that
$ \P B(\partial_\xx)y_1=0$, $\P B(\partial_\xx)z_1 = \P A(\partial_\xx)z_1 $ and hence
 \begin{align}
 \label{PDE.y1}
 \pt y_1
 + \frac{1}{\epsi}\P A(\partial_\xx)z_1
 &=
\sum_{\#J=1}  \P T(v_J),
\\
\label{PDE.z1}
 \pt z_1
 + \frac{\ii}{\epsi^2}\L_1 z_1
 + \frac{1}{\epsi}\Pp B(\partial_\xx)(y_1+z_1)
 &=
 \sum_{\#J=1} \Pp T(v_J).
\end{align}
$\L_1$ is not invertible, but the restricted mapping $ \L_\perp: =\L_1 |_{\Pp \IC^n} : \Pp \IC^n \rightarrow \Pp \IC^n $ 
has an inverse $\L_\perp^{-1}=\Pp \L_\perp^{-1} \Pp:  \Pp \IC^n \rightarrow \Pp \IC^n$.
Therefore,
\eqref{PDE.z1} is equivalent to
 \begin{align}
 \label{PDE.z1.modified}
z_1
 &=
 \ii \epsi \L_\perp^{-1} \Pp \Big(
 \epsi \pt z_1 + B(\partial_\xx)(y_1+z_1) -\epsi \sum_{\#J=1} T(v_J) 
 \Big).
\end{align}

% --------------------------------------------------------------------------------
\subsection{Asymptotic expansion in powers of $\eps$} \label{Ansatz.part.2}
% --------------------------------------------------------------------------------

Our next goal is to approximate a solution of the PDE system \eqref{PDE.vj} with functions $\widetilde v_j$ that are appropriately bounded together with their partial derivatives and that satisfy the initial condition $\P \widetilde v_1(0,\cdot)=\pp$.
The condition \eqref{Ass.initial.data}
suggests an ansatz as a truncated power series in $\eps$ of the form
\begin{align}
\label{Def.vtilde}
\widetilde{v}_j(\tt,\xx) =
\sum_{\ell=j-1}^{m} \epsi^{\ell} v_j^{\ell}(\tt,\xx),
\qquad
j \in \J_+,
\end{align}
where $m$ is the same number as in the definition of $\J$ in \eqref{Def.J}. Note that $\ell$ is a power on $\eps$ and is a superscript on $v_j$.
As before we let $ v_{-j}^{\ell} = \overline{v_j^{\ell}} $.
For example, for $m=5$ we have $\J_+ = \{1, 3, 5 \} $ and \eqref{Def.vtilde} reads
\begin{align}
v_1 &\approx v_1^0 + \epsi v_1^1 + \epsi^2 v_1^2 + \epsi^3 v_1^3 + \epsi^4 v_1^4 + \epsi^5 v_1^5,
\\
\notag
v_3 &\approx \epsi^2 v_3^2 + \epsi^3 v_3^3 + \epsi^4 v_3^4 + \epsi^5 v_3^5,
\\
\notag
v_5 &\approx \epsi^4 v_5^4 + \epsi^5 v_5^5.
\end{align}
As in \eqref{Def.y1.z1} we set
\begin{align}
\label{Def.y1ell.z1ell}
y_1^\ell &= \P v_1^\ell, \quad y_{-1}^\ell = \overline{y_1^\ell}, \qquad 
z_1^\ell = \Pp v_1^\ell, \quad z_{-1}^\ell = \overline{z_1^\ell}
\intertext{and in accordance with \eqref{Ass.initial.data} }
\label{Expansion.y1}
\widetilde y_1 &= \ \P \widetilde v_1 \: = y_1^0 + \epsi y_1^1 + \epsi^2 y_1^2 + \ldots + \epsi^m y_1^m
\\
\label{Expansion.z1}
\widetilde z_1 &= \Pp \widetilde v_1=\hs{1} 0 \hs{1} + \epsi z_1^1 + \epsi^2 z_1^2 + \ldots + \epsi^m z_1^m.
\end{align}
Negative subscripts mean complex conjugation as before. All other terms are set to zero, i.e.
\begin{align}
\label{Other.terms.zero}
z_{\pm 1}^0 = 0, \qquad v_j^\ell = 0 \quad \text{for } \ell < |j|-1.
\end{align}
In particular, this implies $ v_{\pm 1}^0 = y_{\pm 1}^0 $.
For multi-indices $ J=(j_1,j_2,j_3)\in\J^3 $ and
$ L=(\ell_1,\ell_2,\ell_3)\in \IN_0^3 $ we let
\begin{align*}
v_J^L=\big(v_{j_1}^{\ell_1},v_{j_2}^{\ell_2},v_{j_3}^{\ell_3}\big).
\end{align*} 
We will now derive equations for the coefficient functions $y_1^\ell$ for $0\le \ell\le m$, for $z_1^\ell$ for $1\le \ell\le m$, and for
$v_j^{\ell}$ with odd  $3\le j \le m$ for $j-1 \le \ell \le m$.
In order to fulfill the condition $ \P \uu_1(0,\cdot)=\pp $ (cf.~Property~\ref{thm.mfe.property.Puu} in Theorem~\ref{thm.mfe}) we choose the initial data $y_1^0(0,\cdot) = p$ and $ y_1^\ell(0,\cdot) = 0 $ for $ \ell\geq 1$.

\subsection{Equations for $y_1^0$ and $z_1^1$} 
Substituting the expansions \eqref{Def.vtilde}, \eqref{Expansion.y1}, and \eqref{Expansion.z1}
into \eqref{PDE.y1} and \eqref{PDE.z1.modified} and equating like powers of  $\epsi$ yields
 \begin{align}
 \label{y10}
 \pt y_1^0
 + \P A(\partial_\xx)z_1^1
 &=
\sum_{\substack{\#J=1 \\ |L|_1=0}}  \P T(v_J^L),
\\
\label{z11}
z_1^1
 &=
 \ii \L_\perp^{-1} \Pp B(\partial_\xx)y_1^0,
\end{align}
cf.~\cite[Sect.~2.2]{lannes:11}.
Substituting \eqref{z11} into \eqref{y10} yields
 \begin{align}
 \label{v10.temp}
 \pt y_1^0
 + \ii \P A(\partial_\xx) \L_\perp^{-1} \Pp B(\partial_\xx)y_1^0  
 &=
\sum_{\substack{\#J=1 \\ |L|_1=0}} \P T(v_J^L).
\end{align}
With \eqref{Def.B} it follows that
\begin{align*}
\Pp B(\partial_\xx)y_1^0 = \Pp A(\partial_\xx)y_1^0 - (\cg)_\ell \Pp y_1^0 = \Pp A(\partial_\xx)y_1^0,
\end{align*}
because $\Pp y_1^0=0$ according to \eqref{Def.y1ell.z1ell}.
Lemma 2.12 in \cite{lannes:11}
% \myremark{(dort fehlt wohl das $\P$ auf der rechten Seite)}
yields for every smooth $v: \IR^d\to \IC^n$
\begin{align*}
\P A(\partial_\xx) \L_\perp^{-1} \Pp A(\partial_\xx) \P v(\xx) &= -\frac{1}{2} (\nabla_\xx \cdot H\nabla_\xx) \P v(\xx)
%\\
%&= -\frac{1}{2} \sum_{k=1}^d \sum_{\ell=1}^d H_{k\ell} \partial_k \partial_\ell \P v(\xx)
\end{align*}
where 
$H=\nabla^2 \w_1(\kappa)$ is the Hessian of $ \w_1(\cdot)$ at $\kappa$.
Substituting these equations into \eqref{v10.temp} leads to the nonlinear Schr\"odinger equation
\begin{align}
\label{y10.NLS}
 \pt y_1^0
 -\frac{\ii}{2}\nabla_\xx \cdot H\nabla_\xx y_1^0  
 &=
\sum_{\substack{\#J=1 \\ |L|_1=0}} \P T(v_J^L).
\end{align}
Since
 \begin{align*}
\sum_{\substack{\#J=1 \\ |L|_1=0}} T(v_J^L)
&=
T(y_1^0,y_1^0,y_{-1}^0) + T(y_1^0,y_{-1}^0,y_1^0) + T(y_{-1}^0,y_1^0,y_1^0) 
\end{align*}
depends only on $y_1^0$ and $y_{-1}^0=\overline{y_1^0}$, \eqref{y10.NLS} is independent of $ z_1^1$.
Thus, one can first solve \eqref{y10.NLS} with initial value $y_1^0(0,\cdot)=p$
to determine $y_1^0$. This is precisely the initial value problem for the NLS approximation considered in Assumption~\ref{Ass.NLS}.
Therefore, by this assumption a smooth solution $y_1^0=U$ exists over the interval $[0,\tend]$.
Once $y_1^0$ is obtained, $ z_1^1$ is determined by~\eqref{z11}.

\subsection{Equations for $y_1^1$ and $z_1^2$} \label{subsec.y11}

Proceeding in the same way yields
\begin{align}
\label{y11}
 \pt y_1^1 + \P A(\partial_\xx)z_1^2
 &=
 \sum_{\substack{\#J=1 \\ |L|_1=1}}  \P T(v_J^L),
 \\
 \label{z12}
 z_1^2
 &=
 \ii \L_\perp^{-1} \Pp \Big(
 B(\partial_\xx)(y_1^1+z_1^1) - \sum_{\substack{\#J=1 \\ |L|_1=0}} T(v_J^L)
 \Big).
\end{align}
Inserting \eqref{z12} into \eqref{y11} and proceeding as before gives 
\begin{align}
\notag
\pt y_1^1 
 -\frac{\ii}{2}\nabla_\xx \cdot H\nabla_\xx y_1^1  
 + \ii \P A(\partial_\xx)
 \L_\perp^{-1} \Pp \Big(
 B(\partial_\xx)z_1^1 - \sum_{\substack{\#J=1 \\ |L|_1=0}} T(v_J^L)
 \Big) &
 \\
 \label{y11.LS}
=
 \sum_{\substack{\#J=1 \\ |L|_1=1}}  \P T(v_J^L). &
\end{align}
The sum 
\begin{align*}
\sum_{\substack{\#J=1 \\ |L|_1=1}} T(v_J^L)
&=
T(v_1^1,v_1^0,v_{-1}^0) + T(v_1^1,v_{-1}^0,v_1^0) + T(v_{-1}^0,v_1^1,v_1^0) 
\\[-7mm]
&\quad + T(v_1^0,v_1^1,v_{-1}^0) + T(v_1^0,v_{-1}^0,v_1^1) + T(v_{-1}^0,v_1^0,v_1^1) 
\\
&\quad + T(v_1^0,v_1^0,v_{-1}^1) + T(v_1^0,v_{-1}^1,v_1^0) + T(v_{-1}^1,v_1^0,v_1^0) 
\end{align*}
depends only on $y_1^0$ and $v_1^1 = y_1^1 + z_1^1$. Since $v_1^1 $ appears exactly once
in each evaluation of the trilinearity  and since $z_1^1$ does not depend on $ y_1^1$, 
the PDE \eqref{y11.LS} is a \emph{linear} inhomogeneous Schr\"odinger equation for $y_1^1$, to be solved with zero initial value: $y_1^1(0,\cdot)=0$.
With the solution
$y_1^1$ we then compute $z_1^2$ from 
\eqref{z12}.

\subsection{Equations for $y_1^2$, $z_1^3$, and $v_3^2$} 

As in the first two steps one obtains
\begin{align}
 \label{y12}
 \pt y_1^2
 + \P A(\partial_\xx)z_1^3
 &=
\sum_{\substack{\#J=1 \\ |L|_1=2}}  \P T(v_J^L),
\\
\label{z13}
z_1^3
 &=
 \ii \L_\perp^{-1} \Pp \Big(
 \pt z_1^1 + B(\partial_\xx)(y_1^2+z_1^2) - \sum_{\substack{\#J=1 \\ |L|_1=1}} T(v_J^L) 
 \Big).
\end{align}
For the first time there is now a contribution from the time-derivative $\pt z_1$ on the right-hand side of \eqref{PDE.z1.modified}. Substituting \eqref{z13} into \eqref{y12} leads to
 \begin{align}
\notag
 \pt y_1^2
 -\frac{\ii}{2}\nabla \cdot H\nabla y_1^2
 + \ii \P A(\partial_\xx)
  \L_\perp^{-1} \Pp \Big(
 \pt z_1^1 + B(\partial_\xx)z_1^2 - \sum_{\substack{\#J=1 \\ |L|_1=1}} T(v_J^L) 
 \Big) &
 \\
  \label{y12.LS}
 =
\sum_{\substack{\#J=1 \\ |L|_1=2}}  \P T(v_J^L), &
\end{align}
which is to be solved with zero initial value.
In contrast to the first two steps the sum
\begin{align*}
\sum_{\substack{\#J=1 \\ |L|_1=2}} \P T(v_J^L) 
\end{align*}
on the right-hand side involves not only terms $ v_{\pm 1}^\ell $ with subscript $\pm1$, but also 
$ \P T(v_J^L) = \P T(v_3^2,v_{-1}^0,v_{-1}^0) $ and two other terms with permuted arguments. Here,
$ v_3^2 $ has to be computed from \eqref{PDE.vj.modified}. 
Since by construction $ \partial_\tt v_3^0=0 $ and $ v_3^1=0 $ we obtain
\begin{align}
\label{v32}
 v_3^2
 &=
 - \ii \L_3^{-1}\sum_{\substack{\#J=3 \\ |L|_1=0}} T(v_J^L).
 \end{align}
The sum
\begin{align*}
\sum_{\substack{\#J=3 \\ |L|_1=0}} T(v_J^L) = T(v_1^0,v_1^0,v_1^0) = T(y_1^0,y_1^0,y_1^0)
\end{align*}
depends only on $y_1^0$, which is already available.

\subsection{Equations for $y_1^\ell$, $z_1^{\ell+1}$, and $v_j^\ell$} 
For $\ell=3, \ldots, m$ we obtain $y_1^\ell$ by solving the PDE
 \begin{align}
 \notag
 \pt y_1^\ell
 -\frac{\ii}{2}\nabla \cdot H\nabla y_1^\ell
 &+ \P A(\partial_\xx)
 \Big(
  \ii \L_\perp^{-1} \Pp \Big(
 \pt z_1^{\ell-1} + B(\partial_\xx)z_1^\ell - \sum_{\substack{\#J=1 \\ |L|_1=\ell-1}} T(v_J^L) 
 \Big)
 \Big) 
 \\
 \label{y1ell}
 &=
\sum_{\substack{\#J=1 \\ |L|_1=\ell}}  \P T(v_J^L)
\end{align}
with zero initial data. Then, we set
\begin{align}
\label{z1ell}
z_1^{\ell+1}
 &=
 \ii \L_\perp^{-1} \Pp \Big(
 \pt z_1^{\ell-1} + B(\partial_\xx)(y_1^\ell+z_1^\ell) - \sum_{\substack{\#J=1 \\ |L|_1=\ell-1}} T(v_J^L) 
 \Big)
 \intertext{and for $j\in\J_+$ with $3\leq j\leq \ell+1$}
\label{vjell}
 v_j^\ell
 &=
 \ii \L_j^{-1}
 \Big(\pt v_j^{\ell-2} + B(\partial_\xx)v_j^{\ell-1} 
 - \sum_{\substack{\#J=j \\ |L|_1=\ell-2}} T(v_J^L)\Big).
 \end{align}
If $\ell=m$, then \eqref{z1ell} need not be computed, because $ z_1^{m+1} $ appears neither in the ansatz \eqref{Def.vtilde} nor in the construction of the other quantities.

\bigskip

Substituting \eqref{Def.vtilde} into \eqref{Ansatz.v} yields a modulated Fourier expansion
\begin{align}
\label{Approx.mfe}
%\vv(\tt,\xx) \approx 
 \vvtilde(\tt,\xx) &:= 
 \sum_{j\in \J}
e^{\ii j \kappa \cdot \xx/\epsi} 
e^{\ii j (\kappa \cdot c_g - \w) \tt/\epsi^2} \widetilde{v}_j(\tt,\xx)
\\
\label{Approx.mfe.no.vtilde}
&=
\sum_{j\in \J}
e^{\ii j \kappa \cdot \xx/\epsi} 
e^{\ii j (\kappa \cdot c_g - \w) \tt/\epsi^2} 
\sum_{\ell=|j|-1}^{m} \epsi^{\ell} v_j^{\ell}(\tt,\xx),
\end{align}
with $v_{-j}^\ell = \overline{v_j^\ell}$, as an approximation to the solution of \eqref{PDE.vv}.

% --------------------------------------------------------------------------------
\section{Bounds for the modulation functions}
\label{Sec.bounds}
% --------------------------------------------------------------------------------

In this section we derive the bounds \eqref{Ass.initial.data} for the modulation functions $\widetilde v_j$ constructed in the previous section.
A suitable function space for the analysis is the Wiener algebra, which we recall now.

% --------------------------------------------------------------------------------
\subsection{Wiener algebra}
% --------------------------------------------------------------------------------

The Wiener algebra is defined by
\begin{align*}
W &= \{f \in \S'(\IR^d): \widehat{f}\in L^1(\IR^d)\}, \qquad
\|f\|_W = \|\widehat{f}\|_{L^1} = \intl_{\IR^d} |\widehat{f}(k)|_2 \; dk,
\end{align*}
where $\F f = \widehat{f} $ is the Fourier transform of $f$, i.e.
\begin{align*}
(\F f)(k) := (2\pi)^{-d/2}\intl_{\IR^d} f(x)e^{-\ii k\cdot x} dx.
\end{align*}
$W(\IR^d)$ is a Banach algebra and continuously embedded in $L^\infty(\IR^d)$. 
There is a constant $C_T$ such that the trilinear estimate
\begin{align}
\label{Trilinear.estimate}
\| T(f_1,f_2,f_3) \|_W \leq C_T \| f_1 \|_W \| f_2 \|_W \| f_3 \|_W
\end{align}
holds for all $ f_1, f_2, f_3 \in W $; cf.~Section 3.2 in \cite{lannes:11} and p.~715 in \cite{colin-lannes:09}.
 If, in addition, $ g_1, g_2, g_3 \in W $, then
it follows from \eqref{Trilinear.estimate} that
\begin{align}
\label{Trilinear.estimate.difference}
\| T(f_1,f_2,f_3) - T(g_1,g_2,g_3) \|_W 
\leq C_T M^2 \sum_{\ell=1}^3\| f_\ell - g_\ell \|_W
\end{align}
with $M := \max\{ \|f_\ell\|_W, \|g_\ell\|_W, \; \ell=1, 2, 3\}$.
For $s\in\IN$ we define
\begin{align*}
W^s &= \{f \in W(\IR^d): \partial^\alpha f \in W(\IR^d) \text{ for all }
\alpha \in \IN_0^d, |\alpha|_1 \leq s\}, \\
\|f\|_{W^s} &= \sum_{|\alpha|_1 \leq s} \|\partial^\alpha f\|_W.
\end{align*}
In these spaces \eqref{Trilinear.estimate} and the product rule of differentiation imply the trilinear estimates
\begin{align}
\label{Trilinear.estimate.s}
\| T(f_1,f_2,f_3) \|_{W^s} \leq C_T^{(s)} \sum_{s_1+s_2+s_3 = s}\| f_1 \|_{W^{s_1}} \| f_2 \|_{W^{s_2}} \| f_3 \|_{W^{s_3}}
\end{align}
for all $ f_1, f_2, f_3 \in W^s $.

% --------------------------------------------------------------------------------
\subsection{Boundedness of the coefficient functions}
\label{Subsec.bounds}
% --------------------------------------------------------------------------------

The right-hand side of \eqref{Approx.mfe.no.vtilde} can provide a reasonable approximation only if all coefficient functions $ y_1^\ell$, $z_1^\ell$, $v_j^\ell$ remain bounded on $[0,\tend]$. This issue is studied now.
We will work with the regularity
\begin{align} \label{y-reg}
 y_1^0 \in X^s := \bigcap_{i=0}^{\lfloor s/2 \rfloor} C^i([0,\tend],W^{s-2i})
\end{align}
for a sufficiently large integer $s$, where $ \lfloor s/2 \rfloor $ denotes the largest integer which is not larger than $ s/2 $. 
Since  $ y_1^0$ solves the nonlinear Schr\"odinger equation \eqref{y10.NLS}, 
the regularity \eqref{y-reg} with an arbitrary $s$ is implied by Assumption~\ref{Ass.NLS}.

\begin{Lemma}
\label{Lem.everything.is.bounded}
For every integer $q\ge 0$ there is a constant $C$ independent of $\epsi$ such that
\begin{align*}
 \| v_j^{\ell}\|_{X^q} \leq C
\end{align*}
for all $ j\in\J$ and all $ \ell=0, \ldots, m$. 
\end{Lemma}

In particular, this implies the weaker bounds 
$$
 \| v_j^{\ell}\|_{L^\infty([0,\tend] \times \IR^d)} \le 
  \| v_j^{\ell}\|_{L^\infty([0,\tend],W)} \leq   \| v_j^{\ell}\|_{L^\infty([0,\tend],W^q)} \leq C.
 $$

\proof
The first step in the construction of Section~\ref{Ansatz.part.2} is to compute $ y_1^0$ by solving the nonlinear Schr\"odinger equation \eqref{y10.NLS}. 
All other coefficient functions can be traced back to $ y_1^0$, and thus their regularity depends on the regularity of $y_1^0$, as we will explain now. 
Since $ z_1^1 $ is given by \eqref{z11}, it follows from \eqref{y-reg} that $ z_1^1 \in X^{s-1}. $
In Section~\ref{subsec.y11} the next step was to construct $ y_1^1 $ by solving the linear inhomogeneous Schr\"odinger equation \eqref{y11.LS} with zero initial data. The inhomogeneity includes the term
\begin{align}
 - \ii \P A(\partial_\xx) \L_\perp^{-1} \Pp B(\partial_\xx)z_1^1 \in X^{s-3}
\end{align}
and other terms with higher regularity.
By standard regularity theory (e.g. based on Section 4.2~in \cite{pazy:83}), this yields
$y_1^1 \in X^{s-3}$.
Since $z_1^2 $ defined by \eqref{z12} involves $ B(\partial_\xx)y_1^1 $ and other terms of higher regularity, 
we obtain that $z_1^2 \in X^{s-4}.$ 

The coefficient function $ y_1^2 $ is the solution of the linear Schr\"odinger equation 
\eqref{y12.LS} with zero initial data. Here, the critical term in the inhomogeneity is
\begin{align*}
 - \ii \P A(\partial_\xx) \L_\perp^{-1} \Pp  B(\partial_\xx)z_1^2 \in X^{s-6},
\end{align*}
such that $ y_1^2 \in X^{s-6} $. 
An immediate consequence of \eqref{z13} is that $ z_1^3 \in X^{s-7} $.
The function $v_3^2$ defined by \eqref{v32} depends only on $y_1^0 \in X^s$, such that $v_3^2\in W^s.$
By continuing this procedure, it is found that all coefficient functions involved in the approximation
\eqref{Approx.mfe.no.vtilde} are bounded if $s$ is sufficiently large. Roughly speaking, three orders of regularity are required if the upper index of $y_1^\ell$ is increased by one,
i.e. $ y_1^\ell \in X^{s-3\ell}$.
\qed

The bounds for the coefficient functions $v_j^\ell$ in \eqref{Def.vtilde} yield the bounds \eqref{Ass.initial.data} when the truncated expansion $\widetilde v_j$ is taken in the role of $v_j$, as we always do in the following.

% --------------------------------------------------------------------------------
\section{Error bound}
\label{Sec.error}
% --------------------------------------------------------------------------------

We will prove that the modulated Fourier expansion \eqref{Approx.mfe}
approximates a solution of \eqref{PDE.vv} up to $\Ord{\epsi^{m+1}}$.
Let $r$ be the residual of $\vvtilde$, i.e.
\begin{align}
\label{Def.r}
r = \pt \vvtilde + \frac{1}{\epsi}B(\partial_\xx)\vvtilde +\frac{1}{\epsi^2} E\vvtilde 
- T(\vvtilde,\vvtilde,\vvtilde).
\end{align}
Comparing
\begin{align*}
&
\pt \vvtilde + \frac{1}{\epsi}B(\partial_\xx)\vvtilde +\frac{1}{\epsi^2} E\vvtilde 
\\
&=
\sum_{j\in \J} e^{\ii j \kappa \cdot \xx/\epsi} e^{\ii j (\kappa \cdot c_g - \w) \tt/\epsi^2}
\left(
\pt \widetilde{v}_j + \frac{\ii}{\epsi^2}\L_j \widetilde{v}_j + \frac{1}{\epsi}B(\partial_\xx)\widetilde{v}_j
\right)
\end{align*}
with 
\begin{align*}
T(\vvtilde,\vvtilde,\vvtilde) 
%\\
&= 
\sum_{J\in\J^3}  
e^{\ii \#J (\kappa \cdot \xx)/\epsi} 
e^{\ii \#J (\kappa \cdot c_g - \w) \tt/\epsi^2} 
T(\widetilde{v}_{j_1},\widetilde{v}_{j_2},\widetilde{v}_{j_3})
\\
&= 
\sum_{\substack{j \text{ odd} \\ |j|\leq 3m}}
e^{\ii j (\kappa \cdot \xx)/\epsi} 
e^{\ii j (\kappa \cdot c_g - \w) \tt/\epsi^2} 
\sum_{\#J=j} 
T(\widetilde{v}_{j_1},\widetilde{v}_{j_2},\widetilde{v}_{j_3})
\end{align*}
shows that $r$ has the representation
\begin{align}
\label{Representation.r}
r(\tt,\xx) &= 
\sum_{\substack{j \text{ odd} \\ |j|\leq 3m}}
e^{\ii j (\kappa \cdot \xx)/\epsi} e^{\ii j (\kappa \cdot c_g - \w) \tt/\epsi^2} 
r_j(\tt,\xx),
\end{align}
with $ r_{-j} = \overline{r_j} $ and
\begin{align*}
r_j(\tt,\xx) &=
\begin{cases}
\ds
\pt \widetilde{v}_j + \frac{\ii}{\epsi^2}\L_j \widetilde{v}_j + \frac{1}{\epsi}B(\partial_\xx)\widetilde{v}_j
-
\sum_{\#J=j} T(\widetilde{v}_{j_1},\widetilde{v}_{j_2},\widetilde{v}_{j_3})
& \text{if } 1 \leq j \leq m,
\\
\ds
-
\sum_{\#J=j} T(\widetilde{v}_{j_1},\widetilde{v}_{j_2},\widetilde{v}_{j_3})
& \text{if }  m < j \leq 3m.
\end{cases}
\end{align*}
The defects $r_j$ with $j\leq m$ are caused by solving the PDEs \eqref{PDE.vj}, \eqref{PDE.y1}, \eqref{PDE.z1} only \emph{approximately}. For $j>m$ the defects $r_j$ originate from the fact that only the terms with index $ j \in \J = \{\pm 1, \pm 3, \ldots, \pm m \} $ are used in the ansatz \eqref{Ansatz.v}.

\begin{Lemma}\label{Lemma.bound.r}
Let $q\ge 0$ be arbitrary.
There is a constant $C$ independent of~$\eps$ (but dependent on $q$) such that 
\begin{align*}
\sup_{\tt\in[0,\tend]}\| r(\tt,\cdot) \|_{W^q} \leq C \epsi^{m-1}.
\end{align*}
\end{Lemma}
\proof
We will show that
\begin{align*}
\sup_{\tt\in[0,\tend]}\| r_j(\tt,\cdot) \|_{W^q} \leq C \epsi^{m-1}
\end{align*}
for all odd $ j = 1, \ldots , 3m$ for some constant $C$. Then, the assertion follows directly from
\eqref{Representation.r}.

Lemma~\ref{Lem.everything.is.bounded} and \eqref{Def.vtilde} imply that 
$ \| \widetilde{v}_j \|_{W^q} = \Ord{\epsi^{j-1}}$ and, via the trilinear estimate 
\eqref{Trilinear.estimate.s}, that
\begin{align}
\label{Lemma.bound.r.Eq01}
\| T(\widetilde{v}_{j_1},\widetilde{v}_{j_2},\widetilde{v}_{j_3}) \|_{W^q}
\leq c \, \epsi^{j_1-1}\epsi^{j_2-1}\epsi^{j_3-1} = c \, \epsi^{\#J-3}
\end{align}
for some constant $c$ which depends on $C_T^{(q)}$ from \eqref{Trilinear.estimate.s}.
If $ j>m $, then $ j\geq m+2$ because $j$ and $m$ are both odd numbers, and we obtain
\begin{align*}
\| r_j(t,\cdot) \|_{W^q} \leq \sum_{\#J=j} 
\| T(\widetilde{v}_{j_1},\widetilde{v}_{j_2},\widetilde{v}_{j_3}) \|_{W^q}
\leq C\epsi^{j-3} \leq C \epsi^{m-1},
\end{align*}
where $C$ depends on $c$ from \eqref{Lemma.bound.r.Eq01}
and on the number of multi-indices $J\in\J^3$ with $\#J=j$.

Now let $ j \leq m$. In this case, $r_j$ is given by
\begin{align}
\label{Residual.expansion}
r_j(\tt,\xx) &= \sum_{\ell=j-1}^m \epsi^\ell r_j^\ell(\tt,\xx)
\intertext{with}
\notag
r_j^\ell(\tt,\xx) &= 
\pt v_j^\ell + \ii\L_j v_j^{\ell+2} + B(\partial_\xx)v_j^{\ell+1}
-
\sum_{\substack{\#J=j \\ |L|_1=\ell}}  T(v_J^L).
\end{align}
In Subsection~\ref{Ansatz.part.2} the coefficients $ y_1^\ell$, $z_1^\ell$, $v_j^\ell$ have been constructed in such a way that $r_j^\ell = 0 $ for $ \ell \leq m-2 $. Hence, only
\begin{subequations}
\label{Lemma.bound.r.Eq02}
\begin{align}
r_j^{m-1} &= 
\pt v_j^{m-1} + B(\partial_\xx)v_j^m - \sum_{\substack{\#J=j \\ |L|_1=m-1}}  T(v_J^L),
\\
r_j^m &=
\pt v_j^m - \sum_{\substack{\#J=j \\ |L|_1=m}}  T(v_J^L)
\end{align}
\end{subequations}
do not vanish. 
Since $ \| r_j^{m-1} \|_{W^q} $ and $ \| r_j^{m} \|_{W^q} $ remain uniformly bounded by Lemma~\ref{Lem.everything.is.bounded}, it follows from \eqref{Residual.expansion}
that 
\begin{align*}
\sup_{\tt\in[0,\tend]}\| r_j(\tt,\cdot) \|_{W^q} \leq C \epsi^{m-1}
\end{align*}
for all $ m < j \leq 3m $. This proves the assertion.
\qed

\begin{Lemma}\label{Lem.error.bound}
There exists a unique mild solution $\vv$ of \eqref{PDE.vv} with initial data $\vv(0,\cdot)=\widetilde \vv(0,\cdot)$ in the space 
$C(I,W^q)$
(with arbitrary $q\ge 0$) 
on the whole time interval $I=[0,\tend]$ of Assumption~\ref{Ass.NLS}, and 
\begin{align}\label{error.bound}
 \| \vv - \vvtilde \|_{L^\infty([0,\tend],W^q)}
 &\leq C \epsi^{m+1},
\end{align}
where $C$ is independent of $\eps$ (but depends on $m$ and $q$).
\end{Lemma}

In particular, this implies
$$
 \| \vv - \vvtilde \|_{L^\infty([0,\tend]\times \IR^d)} \leq C \epsi^{m+1}.
$$

\proof We split the proof into three parts (a)--(c).

(a) For every $\epsi>0$ the operator
\begin{align*}
\A_\epsi = - \frac{1}{\epsi}B(\partial_\xx) - \frac{1}{\epsi^2} E 
\quad \text{with domain} \quad
D(\A_\epsi) = W^{q+1}
\end{align*}
generates a strongly continuous group $ \big(\exp(\tt \A_\epsi)\big)_{\tt\in\IR} $ on $W^q$.
The group operators are explicitly given by
\begin{align*}
\F\big(\exp(\tt \A_\epsi)f\big)(k) = 
\exp\left( - \frac{\tau}{\epsi^2}\big(\ii \epsi B(k) + E \big)\right)\widehat{f}(k)
\end{align*}
for every $f\in W^q$ (with arbitrary $q\ge 0$) and all $ \tt\in\IR$.
Since $E\in\IR^{n\times n}$ is skew-symmetric and $ \ii B(k) $ is skew-Hermitian, it follows that
\begin{align*}
\| \exp(\tau \A_\epsi)f \|_{W^q} 
&= \sum_{|\alpha|_1\leq q}\big\| \F\big(\exp(\tau \A_\epsi)\partial_\xi^\alpha f\big) \big\|_{L^1}
\\
&= \sum_{|\alpha|_1\leq q} \
\intl_{\IR^d} \left|\exp\Big(-\tfrac{\tau}{\epsi^2}\big(\ii \epsi B(k) + E\big) \Big)k^\alpha \widehat{f}(k)\right|_2 \; dk
\\
&= \sum_{|\alpha|_1\leq q} \
\intl_{\IR^d} |k^\alpha \widehat{f}(k)|_2 \; dk
\\
&= \sum_{|\alpha|_1\leq q} \| \partial_\xi^\alpha f \|_{W^q}  = \| f \|_{W^q}
\end{align*}
so that $ \exp(t \A_\epsi): W^q \rightarrow W^q $ is an isometry.
Together with the trilinear estimate \eqref{Trilinear.estimate.s}, this yields that a 
unique mild solution $\vv$ of \eqref{PDE.vv} with initial data $\vv(0,\cdot)=\widetilde \vv(0,\cdot)\in W^q$ exists in the space 
$C([0,\bar\tau],W^q)$ for every $\overline\tau < \tau_{\max}$, for some maximal time $\tau_{\max}$ where the solution becomes unbounded. We will find later that 
$\tau_{\max}>\tend$.

(b) By Lemma~\ref{Lem.everything.is.bounded}, $ \| \vvtilde \|_{L^\infty([0,\tend],W^q)}
\le C$. 
Comparing \eqref{Def.r} with \eqref{PDE.vv} shows that the error
$ \err=\vv-\vvtilde $ solves the evolution equation
\begin{align}
 \label{PDE.delta}
\partial_\tt \err =  
- \frac{1}{\epsi}B(\partial_\xx)\err
- \frac{1}{\epsi^2} E\err
+ \left[T(\vv,\vv,\vv)-T(\vvtilde,\vvtilde,\vvtilde)\right]
- r.
\end{align}

Applying Duhamel's principle (variation of constants formula) to \eqref{PDE.delta} and using the trilinear estimate \eqref{Trilinear.estimate.difference} yields
\begin{align*}
\| \err(\tt) \|_{W^q}  
&\leq \| \err(0) \|_{W^q}  
+ \widetilde C \intl_0^\tt \|\err(\sigma)\|_{W^q}  \; d\sigma
+ \intl_0^\tt \| r(\sigma) \|_{W^q}  \; d\sigma.
\end{align*}
Since $\err(0)=0$ by assumption and since $ \| r(\sigma) \|_{W^q}  \leq C\epsi^{m-1} $ by Lemma~\ref{Lemma.bound.r}, we obtain
\begin{align*}
\| \err(\tt) \|_{W^q}  
&\leq
\widetilde C  \intl_0^\tt \|\err(\sigma)\|_{W^q}  \; d\sigma
+ C\tend \epsi^{m-1},
\end{align*}
such that the error bound \eqref{error.bound} follows from Gronwall's lemma. This holds true on subintervals of $[0,\tend]$ as long as the solution $v\in C([0,\overline \tau], W^q)$ exists, i.e. for $\overline\tau < \tau_{\max}$. Since $\| v(\tau) \|_{W^q} \le \|v- \widetilde v \|_{W^q} + \| \widetilde v \|_{W^q}$,
the solution remains bounded and hence can be continued up to and beyond $\overline\tau=\tend$, so that $\tau_{\max}>\tend$.

\medskip
(c) So far, we have only obtained an $\Ord{\eps^{m-1}}$ error bound instead of the stated $\Ord{\eps^{m+1}}$ bound. The improvement of the order of approximation comes about as follows: We can repeat the same argument with $m+2$ instead of~$m$. The resulting error,
denoted $e_{m+2}$, is then of size $\Ord{\eps^{m+1}}$ and differs from the above error $e_m=e$ (for the expansion up to terms of order $\Ord{\eps^m}$) only by the
terms $\eps^\ell v_j^{\ell}$ with $\ell=m+1,m+2$, i.e.
\begin{align*}
e_m(\tt,\xx) &= e_{m+2}(\tt,\xx) 
+ \sum_{|j|\le m+2} 
e^{\ii j \kappa \cdot \xx/\epsi} 
e^{\ii j (\kappa \cdot c_g - \w) \tt/\epsi^2} 
\hs{-5} \sum_{\ell\in\{m+1,m+2\}} \hs{-5}\epsi^{\ell} v_j^{\ell}(\tt,\xx).
\end{align*}
By Lemma~\ref{Lem.everything.is.bounded} and the bound for $e_{m+2}$, each of these terms is $\Ord{\eps^{m+1}}$, and so we obtain the stated error bound.
\qed

Finally we have proved \eqref{Ansatz.v}--\eqref{v-err}. By returning to the original coordinates $(t,x)$, we thus obtain items 1.--\,4. of Theorem~\ref{thm.mfe}. The uniqueness of the modulation functions $u_j$ up to $O(\eps^{m+1})$ follows from the uniqueness of the functions $v_j^\ell$ in the construction of Section~\ref{Sec.construction}. This completes the proof of Theorem~\ref{thm.mfe}.

\bibliographystyle{alpha}
\bibliography{literature}

\end{document}